\input amstex \loadmsbm
\documentstyle{amsppt}
\magnification=1200 \nopagenumbers \baselineskip=13pt
\parskip=4pt plus 1pt
\topskip=35pt
\hsize 15.5truecm
\vsize 22.5truecm

\font\ver=cmcsc10

\font\aut=cmr12

 \font\rms=cmr8 \font\rm=cmr9 \font\bfs=cmbx8 \font\bfm=cmbx9 \font\bfg=cmbx10 \font\its=cmti8
 \font\tit=cmss17


\def\headimpar{\hfil\ver Purely infinite simple Leavitt path algebras
\hfil{\rm\folio}}
\def\headpar{{\rm \folio}\hfil\ver Abrams and Aranda \hfil}
\headline={\ifnum\pageno=1\hfil
\else\ifodd\pageno\headimpar\else\headpar\fi\fi}
\footline={\ifnum\pageno=1\hfil\tenrm{\rm \folio}\hfil\else\hfil\fi}


\newcount\parracount
\newcount\recount
\parracount=0
\recount=0



\def\re#1#2{\advance\recount by 1 \edef#1{\number\recount}
\def\nnnn{#2} \ifx\nnnn\empty\vskip 3 pt plus 1 pt{\bf #1. }\else \vskip 3 pt plus 1
pt{ \ver #2 }{\bf #1.}\it\fi }

\def\redef#1#2{\advance\recount by 1 \edef#1{\number\recount}
\def\nnnn{#2} \ifx\nnnn\empty\vskip 3 pt plus 1 pt{\bf #1. }\else \vskip 3 pt plus 1
pt{ \ver #2 }{\bf #1.}\fi }

\def\endre{\par\rm}
\def\endredef{\par\rm}

\def\proof{{\bf Proof:} }
\def\endproof{~\vrule height 1.2ex width 1.2ex depth -0.06ex \par}
\def\parrafo#1{\advance\parracount by 1\recount=0
\vfil\vskip 10 pt plus 10 pt\noindent{\bf \number\parracount. #1}\par\nobreak\vskip
3 pt  plus 3 pt\nobreak}





\def\AA{1}
\def\AM{2}
\def\Ados{3}
\def\AGGP{4}
\def\AGP{5}
\def\AMP{6}
\def\BPRS{7}
\def\KPR{8}
\def\Lone{9}
\def\Ltwo{10}
\def\P{11}


\vskip 20 pt \centerline{\tit Purely infinite simple Leavitt path algebras}

\vskip 14 pt \centerline{\aut Gene Abrams $^{a,}$\plainfootnote{\rms $^*$}{\rms
Corresponding author.} \aut and Gonzalo
Aranda Pino $^b$} \vskip 9 pt \centerline{\its $^a$Department of Mathematics,
University of Colorado, Colorado Springs,
CO 80933, U.S.A.} \vskip -1 pt \centerline{\its $^b$Departamento de \'Algebra,
Geometr\'\i a y Topolog\'\i a,
Universidad de M\'alaga, 29071 M\'alaga, Spain} \plainfootnote{}{{\its E-mail
addresses:} {\rms abrams\@math.uccs.edu
(G. Abrams), gonzalo\@agt.cie.uma.es (G. Aranda).}}


\vskip 15 pt

{\narrower\noindent\baselineskip=5 pt {\bfs Abstract}

\rms We give necessary and sufficient conditions on a row-finite graph E so that the Leavitt path algebra L(E) is
purely infinite simple.  This result provides the algebraic analog to the corresponding result for the Cuntz-Krieger
C$^*$-algebra C$^*$(E) given in \cite{\BPRS}.}

\medskip {{\its Keywords:} {\rms purely infinite; path algebra; Leavitt algebra;
Cuntz Krieger C*-algebra.}}




\endre

An idempotent $e$ in a ring $R$ is called {\bfm infinite} if $eR$ is isomorphic as a right $R$-module to a proper
direct summand of itself.  $R$ is called {\bfm purely infinite} in case every right ideal of $R$ contains an infinite
idempotent.   Much recent attention has been paid to the structure of purely infinite simple rings, from both an
algebraic (see e.g. \cite{\Ados}, \cite{\AGGP}, \cite{\AGP}) as well as an analytic (see e.g. \cite{\BPRS},
\cite{\KPR}, \cite{\P}) point of view.  The Leavitt path algebra $L(E)$ of a graph $E$ is investigated in \cite{\AA}.
$L(E)$ is the algebraic counterpart of the Cuntz-Krieger algebra $C^*(E)$; furthermore, the class of algebras of the
form $L(E)$ significantly broadens the collection of algebras studied by Leavitt in his seminal papers \cite{\Lone} and
\cite{\Ltwo}.    In \cite{\AA} the authors give necessary and sufficient conditions on $E$ so that $L(E)$ is simple. In
the current article we provide necessary and sufficient conditions on $E$ so that $L(E)$ is purely infinite simple
(Theorem 11).

We recall the definition of the Leavitt path algebra $L(E)$.

\redef{\pathalgebra}{\bf Definitions} A {\bfm (directed) graph} $E=(E^0,E^1,r,s)$ consists of two countable sets
$E^0,E^1$ and functions $r,s:E^1 \to E^0$. The elements of $E^0$ are called {\bfm vertices} and the elements of $E^1$
{\bfm edges}. Let $K$ be a field.  The {\bfm path $K$-algebra over $E$} is the free associative $K$-algebra $K[E^0\cup
E^1]$ with relations given by: $v_iv_j=\delta_{ij}v_i$ for every $v_i,v_j\in E^0$, and $e_i=e_ir(e_i)=s(e_i)e_i$ for
every $e_i\in E^1.$ The {\bfm extended graph of} $E$ is the graph $\widehat{E}=(E^0,E^1\cup (E^1)^*,r',s')$, where
$(E^1)^*=\{e_i^*:e_i\in E^1\}$ and the functions $r'$ and $s'$ are defined as: $r'|_{E^1}=r,\ s'|_{E^1}=s,\
r'(e_i^*)=s(e_i)$, and $s'(e_i^*)=r(e_i).$ We call the elements of $E^1$ (resp., $(E^1)^*$) the {\bfm real edges}
(resp., the {\bfm ghost edges}) of $E$.

Now suppose that $E$ is {\bfm row-finite} (i.e., that $s^{-1}(v)$ is finite for all $v \in E^0$.) The {\bfm Leavitt
path algebra of $E$ with coefficients in $K$}, denoted by $L_K(E)$ (or $L(E)$ when appropriate), is defined as the path
$K$-algebra over the extended graph $\widehat{E}$, satisfying the so-called {\bfm Cuntz-Krieger relations}:
\item{} (CK1) $e_i^*e_j=\delta_{ij}r(e_j)$ for every $e_j\in E^1$ and  $e_i^*\in
(E^1)^*$, and \item{} (CK2) $v_i=\sum_{\{e_j\in E^1:s(e_j)=v_i\}}e_je_j^*$ for every  $v_i\in E^0$ for which
$s^{-1}(v_i) \neq \emptyset$.
\endredef

\medskip

\re{\Leavittalgebraexamples}{\bf Examples}

{\rm \item{(i)} Let $E$ be the ``finite line" graph defined by $E^0=\{v_1,\dots,v_n\}$, $E^1=\{y_1,\dots,y_{n-1}\}$,
$s(y_i)=v_i$, and $r(y_i)=v_{i+1}$ for $i=1,\dots,n-1$.  Then $L(E) \cong M_n(K)$, via the map $v_i \mapsto e_{ii}, y_i
\mapsto e_{ii+1}$, and $y_i^* \mapsto e_{i+1i}$ (where $e_{ij}$ denotes the standard $(i,j)$-matrix unit in $M_n(K)$).
\item{(ii)} Let $n\geq 2$. Let $E$ be the ``rose with $n$ leaves" graph defined by
$E^0=\{*\}$, $E^1=\{y_1,\dots,y_n\}$. Then $L(E)\cong L(1,n)$, the {\it Leavitt algebra} investigated in \cite{\Ltwo}.
Specifically, $L(E)$ is isomorphic to the free associative $K$-algebra with generators $\{x_i,y_i:1\leq i\leq n\}$ and
relations
$$ (1)\  x_iy_j=\delta_{ij} \hbox{ for all } 1\leq i,j\leq n, \hbox{ and } (2)\ \sum_{i=1}^n y_ix_i=1 .$$
\endre}
\endre

Throughout this article all graphs will be assumed to be row-finite. We briefly establish some graph-theoretic
notation.  For each edge $e$, $s(e)$ is the {\bfm source} of $e$ and $r(e)$ is the {\bfm range} of $e$. A vertex $v$
for which $s^{-1}(v) = \emptyset$ is called a {\bfm sink}. A graph $E$ is {\bfm finite} if $E^0$ is a finite set.
 A {\bfm path} $\mu$ in a graph $E$ is a sequence of edges $\mu=\mu_1 \dots
\mu_n$ such that $r(\mu_i)=s(\mu_{i+1})$ for $i=1,\dots,n-1$. In such a case, $s(\mu):=s(\mu_1)$ is the source of $\mu$
and $r(\mu):=r(\mu_n)$ is the range of $\mu$. For vertices we define $r(v)=v=s(v)$.  We define a preorder $\leq$ on
$E^0$ given by: $v\leq w$ in case $w=v$ or there is a path $\mu$ such that $s(\mu)=v$ and $r(\mu)=w.$ If
$s(\mu)=r(\mu)$ and $s(\mu_i)\neq s(\mu_j)$ for every $i\neq j$, then $\mu$  is a called a {\bfm cycle}. $E$ is {\bfm
acyclic} if $E$ contains no cycles.  The set of paths of length $n>0$ is denoted by $E^n$. The set of all paths (and
vertices) is $E^*:=\cup_{n\geq 0} E^n$. It is shown in \cite{\AA} that $L(E)$ is a ${\Bbb Z}$-graded $K$-algebra,
spanned as a $K$-vector space by $\{pq^* \mid p,q$ are paths in $E \}$. By \cite{\AA, Lemma 1.6}, $L(E)$ is unital if
and only if $E$ is finite; otherwise, $L(E)$ is a ring with set of local units consisting of sums of distinct vertices.

\medskip

If $\alpha \in L(E)$ and $d \in {\Bbb Z}^+$, then we say that $\alpha$ is {\bfm
representable as an element of degree $d$ in real (resp. ghost) edges}
in case $\alpha$ can be written as a sum of monomials from the aforementioned
spanning set of $L(E)$, in such a way that $d$ is the maximum length of a path $p$
(resp. $q$) which appears in such monomials.   We note that an element of $L(E)$ may
be representable as an element of different degrees in real (resp. ghost) edges,
depending on the particular representation used for $\alpha$.

\re{\finitedimensional}{\bf Lemma} Let $E$ be a finite acyclic graph.  Then $L(E)$ is finite dimensional.
\endre
\proof Since the graph is row-finite, the given condition on $E$ is equivalent to the condition that $E^*$ is finite.
 The result now follows from the previous observation
that $L(E)$ is spanned as a $K$-vector space by $\{pq^* \mid p,q$ are paths in $E \}$.
\endproof

Lemma \finitedimensional\ is precisely the tool we need to establish the following key result.

\re{\ultramatricial}{\bf Proposition} Let $E$ be a graph.  Then $E$ is acyclic if and only if $L(E)$ is a union of a
chain of finite dimensional subalgebras.
\endre
\proof Assume first that $E$ is acyclic.   If $E$ is finite, then Lemma 3 gives the result. So now suppose $E$ is
infinite, and rename the vertices of $E^0$ as a sequence $\{v_i\}_{i=1}^\infty$. We now define a sequence
$\{F_i\}_{i=1}^\infty$ of subgraphs of $E$. Let $F_i=(F_i^0,F_i^1,r,s)$ where $F_i^0:=\{v_1,\dots,v_i\}\cup
r(s^{-1}(\{v_1,\dots,v_i\})$, $F_i^1:=s^{-1}(\{v_1,\dots,v_i\})$, and $r,s$ are induced from $E$.     In particular,
$F_i \subseteq F_{i+1}$ for all $i$. For any $i>0$, $L(F_i)$ is a subalgebra of $L(E)$ as follows. First note that we
can construct $\phi:L(F_i) \to L(E)$ a $K$-algebra homomorphism because the Cuntz-Krieger relations in $L(F_i)$ are
consistent with those in $L(E)$, in the following way: Consider $v$ a sink in $F_i$ (which need not be a sink in $E$),
then we do not have CK2 at $v$ in $L(F_i)$. If $v$ is not a sink in $F_i$, then there exists $e\in
F_i^1:=s^{-1}(\{v_1,\dots,v_i\})$ such that $s(e)=v$. But $s(e)\in \{v_1,\dots,v_i\}$ and therefore $v=v_j$ for some
$j$, and then $F_i^1:=s^{-1}(\{v_1,\dots,v_i\})$ ensures that all the edges coming to $v$ are in $F_i$, so CK2 at $v$
is the same in $L(F_i)$ as in $L(E)$. The other relations offer no difficulty.  Now, with a similar construction and
argument to that used in \cite{\AA, Proof of Theorem 3.11} we find $\psi:L(E) \to L(F_i)$ a $K$-algebra homomorphism
such that $\psi \phi=Id|_{L(F_i)}$, so that $\phi$ is a monomorphism, which we view as the inclusion map. By
construction, each vertex in $E^0$ is in $F_i$ for some $i$; furthermore, the edge $e$ has $e\in F_j^1$, where
$s(e)=v_j$. Thus we conclude that $L(E)= \cup_{i=1}^\infty L(F_i)$. (We note here that the embedding of graphs
$j:F_i\hookrightarrow E$ is a complete graph homomorphism in the sense of \cite{\AMP}, so that the conclusion $L(E)=
\cup_{i=1}^\infty L(F_i)$ can also be achieved by invoking \cite{\AMP, Lemma 2.1}.)

Since $E$ is acyclic, so is each $F_i$. Moreover, each $F_i$ is finite since, by the row-finiteness of $E$, in each
step we add only finitely many vertices. Thus, by Lemma 3, $L(F_i)$ is finite dimensional, so that $L(E)$ is indeed a
union of a chain of finite dimensional subalgebras.

For the converse, let $p\in E^*$ be a cycle in $E$. Then $\{p^m\}_{m=1}^\infty$ is a linearly independent infinite set,
so that $p$ is not contained in any finite dimensional subalgebra of $L(E)$.
\endproof

We note that when $E$ is finite and acyclic then $L(E)$ can be shown to be isomorphic to a finite direct sum of full
matrix rings over $K$, and, for any acyclic $E$, $L(E)$ is a direct limit of subalgebras of this form.  The proof
follows along the same lines as that given in \cite{\KPR, Corollary 2.2 and 2.3}.

The description of the simple Leavitt path algebras given in \cite{\AA} will play a key role here, so we briefly review
the germane ideas.  An edge $e \in E^1$ is an {\bfm exit} to the path $\mu=\mu_1 ... \mu_n$ if there exists $i$ such
that $s(e)=s(\mu_i)$ and $e\neq \mu_i$.
 A vertex $w\in E^0$ {\bfm connects to} $v\in E^0$ if
$w\leq v$.  A subset $H\subseteq E^0$ is {\bfm hereditary} if $w\in H$ and $w\leq v$ imply $v\in H$; $H$ is {\bfm
saturated} if whenever $s^{-1}(v)\neq \emptyset$ and $\{r(e):s(e)=v\}\subseteq H$, then $v\in H$.  The main result of
\cite{\AA} is the following

\re{\simple}{\bf Theorem} \cite{\AA, Theorem 3.11} Let $E$ be a graph. Then $L(E)$ is simple if and only if:
\item{(i)} The only hereditary and saturated subsets of $E^0$ are $\emptyset$
and $E^0$, and
\item{(ii)} Every cycle in $E$ has an exit.
\endre

The following Proposition is a useful rearrangement of one of the consequences of the proof of Theorem \simple.

\re{\reducetovertex}{\bf Proposition} Let $E$ be a graph with the property that every cycle has an exit. Then for every
nonzero $\alpha\in L(E)$ there exist $a,b\in L(E)$ such that $a\alpha b \in E^0$.
\endre
\proof Let $\alpha$ be representable by an element having degree $d$ in real edges. If $d=0$, then by \cite{\AA,
Corollary 3.7} we are done. So suppose $d>0$. By \cite{\AA, Lemma 1.5}, given a monomial which is not a vertex, either
it begins with a real edge or all its edges are ghost edges. Then we can write $$\alpha =\sum_{n=1}^m e_{i_n}
\alpha_{e_{i_n}} + \beta$$ where $m\geq 1$, $e_{i_n}\alpha_{e_{i_n}}\neq 0$ for every $n$, each $\alpha_{e_{i_n}}$ is
representable as an element of degree less than that of $\alpha$ in real edges, and $\beta$ is a polynomial in only
ghost edges (possibly zero). We will present a process by which we will find $\widehat{a},\widehat{b}$ such that
$\widehat{a}\alpha \widehat{b}\neq 0$ and is representable as an element having degree less than $d$ in real edges.

For an arbitrary edge $e_j\in E^1$, we have two cases:

Case 1: $j\in \{i_1,\dots, i_m\}$. Then $e_j^*\alpha=\alpha_{e_j}+e_j^*\beta$. If
this element is nonzero then by
choosing $\widehat{a}=e_j^*$ and $\widehat{b}$ a local unit for $\alpha$ we would
be done. For later use, we note that if $e_j^*\alpha$ is zero, then
$\alpha_{e_j}=-e_j^*\beta$, and therefore $e_j\alpha_{e_j}=-e_je_j^*\beta$.

Case 2: $j\not\in \{i_1,\dots, i_m\}$. Then $e_j^*\alpha=e_j^*\beta$.  If
$e_j^*\beta \neq 0$, then with
$\widehat{b}$ as before we would have $e_j^*\alpha \widehat{b}$ is a nonzero
polynomial which is representable as an element having degree $0<d$ in real edges,
and again we would be done. For later use, we note that if $e_j^*\beta = 0$, then in particular we
have
$0=-e_je_j^*\beta$.

So we may assume that we are in the latter possibilities of both Case 1 and 2; i.e., we may assume that $e^*\alpha = 0$
for all $e \in E^1$.   We show that this situation cannot happen. First, suppose $v$ is a sink in $E$. Then we may
assume $v\beta =0$, as follows. Multiplying the displayed equation by $v$ on the left gives $v\alpha = v\sum_{n=1}^m
e_{i_n} \alpha_{e_{i_n}} + v\beta$.  Since $v$ is a sink we have $v e_{i_n} =0$ for all $1 \leq n \leq m$, so that
$v\alpha = v\beta $. But if $v\beta \neq 0$ then $\widehat{a}=v$ and $\widehat{b}$ as above would yield a nonzero
element in only ghost edges and we would be done as in Case 2.

Now let $S_1=\{v_j \in E^0:v_j = s(e_{i_n}) \hbox{ for some } 1\leq n \leq m\}$, and let $S_2 =
\{v_{k_1},...,v_{k_t}\}$ where $(\sum_{i=1}^t v_{k_i})\beta = \beta$. We note that $w\beta = 0$ for every $w\in E^0 -
S_2$.  Also, by definition there are no sinks in $S_1$, and by a previous observation we may assume that there are no
sinks in $S_2$.     Let $S = S_1 \cup S_2$. Then in particular we have $(\sum _{v\in S} v)\beta = \beta$.

We now argue that in this situation $\alpha$ must be zero.  To this end, $$\eqalign{
 \alpha & =\sum_{n=1}^m e_{i_n} \alpha_{e_{i_n}} + \beta = \sum_{n=1}^m
-e_{i_n}e_{i_n}^*\beta + \beta \hskip5mm
 \hbox{ (by Case 1)} \cr
 & = \sum_{n=1}^m -e_{i_n}e_{i_n}^*\beta - \left(\sum_{_{\ \ \ s(e_j)\in S}^{j\notin
\{i_1,...,i_m\}}} e_je_j^*\right) \beta  + \beta  \cr
 & \hskip2cm \hbox{ (by Case 2, the newly subtracted terms equal 0) }\cr
 & =-(\sum _{v\in S} v)\beta + \beta  \hskip5mm \hbox{ (no sinks in } S \hbox{
implies that CK2 applies at each }v\in S)
\cr
 & = -\beta + \beta = 0.\cr}
 $$

As we have assumed $\alpha \neq 0$ we have reached the desired contradiction. Thus
we are always able to find
$\widehat{a},\widehat{b}$ such that $\widehat{a}\alpha \widehat{b}$ is nonzero, and
is representable in degree
less than $d$ in real edges.   By repeating this process enough times ($d$ at most),
we can
find $\widehat{a_k}\dots \widehat{a_1}, \widehat{b_1}\dots \widehat{b_k}$ such that
we can represent
$\widehat{a_k}\dots \widehat{a_1}\alpha \widehat{b_1}\dots \widehat{b_k}\neq 0$ by
an element of degree zero in real
edges.  Thus \cite{\AA, Corollary 3.7} applies, and finishes the proof.
\endproof

A {\bfm closed simple path based at $v_{i_0}$} is a path $\mu=\mu_1 \dots \mu_n$, with $\mu_j \in E^1$, $n\geq 1$ such
that $s(\mu_j)\neq v_{i_0}$ for every $j>1$ and $s(\mu)=r(\mu)=v_{i_0}$. Denote by $CSP(v_{i_0})$ the set of all such
paths.  We note that a cycle is a closed  simple path based at any of its vertices, but not every closed simple path
based at $v_{i_0}$ is a cycle.  We define the following subsets of $E^0$:$$\eqalign{
 V_0 & = \{v\in E^0: CSP(v)=\emptyset\}  \cr
 V_1 & = \{v\in E^0: |CSP(v)|=1\}  \cr
 V_2 & = E^0-(V_0 \cup V_1)\cr} $$

\re{\voneempty}{\bf Lemma} Let $E$ be a graph. If $L(E)$ is simple, then
$V_1=\emptyset$.
\endre
\proof For any subset $X\subseteq E^0$ we define the following subsets. $H(X)$ is the set of all vertices that can be
obtained by one application of the hereditary condition at any of the vertices of $X$; that is, $H(X):=r (s^{-1}(X))$.
Similarly, $S(X)$ is the set of all vertices obtained by applying the saturated condition among elements of $X$, that
is, $S(X):=\{v\in E^0: \emptyset \neq \{r(e):s(e)=v\}\subseteq X\}$.  We now define $G_0:=X$, and for $n\geq 0$ we
define inductively $G_{n+1}:=H(G_n)\cup S(G_n) \cup G_n$. It is not difficult to show that the smallest hereditary and
saturated subset of $E^0$ containing $X$ is the set $G(X):=\bigcup_{n\geq 0} G_n$.

Suppose now that $v\in V_1$, so that $CSP(v)=\{p\}$. In this case $p$ is clearly a cycle. By Theorem \simple\ we can
find an edge $e$ which is an exit for $p$. Let $A$ be the set of all vertices in the cycle. Since $p$ is the only cycle
based at $v$, and $e$ is an exit for $p$, we conclude that $r(e)\not\in A$. Consider then the set $X=\{r(e)\}$, and
construct $G(X)$ as described above.  Then $G(X)$ is nonempty and, by construction, hereditary and saturated.

Now Theorem \simple\ implies that $G(X)=E^0$, so we can find $n=min\{m: A\cap G_m \neq \emptyset\}$. Take $w\in A\cap
G_n$. We are going to show that $w\geq r(e)$. First, since $r(e)\not\in A$, then $n>0$ and therefore $w\in
H(G_{n-1})\cup S(G_{n-1}) \cup G_{n-1}$. Here, $w\in G_{n-1}$ cannot happen by the minimality of $n$. If $w\in
S(G_{n-1})$ then $\emptyset\neq \{r(e):s(e)=w\}\subseteq G_{n-1}$. Since $w$ is in the cycle $p$, there exists $f\in
E^1$ such that $r(f)\in A$ and $s(f)=w$. In that case $r(f)\in A\cup G_{n-1}$ again contradicts the minimality of $n$.
So the only possibility is $w\in H(G_{n-1})$, which means that there exists $e_{i_1}\in E^1$ such that $r(e_{i_1})=w$
and $s(e_{i_1})\in G_{n-1}$.

We now repeat the process with the vertex $w'=s(e_{i_1})$. If $w'\in G_{n-2}$ then we would have $w\in G_{n-1}$, again
contradicting the minimality of $n$. If $w'\in S(G_{n-2})$ then, as above, $\{r(e):s(e)=w'\}\subseteq G_{n-2}$, so in
particular would give $w=r(e_{i_1})\in G_{n-2}$, which is absurd. So therefore $w'\in H(G_{n-2})$ and we can find
$e_{i_2}\in E^1$ such that $r(e_{i_2})=w'$ and $s(e_{i_2})\in G_{n-2}$.

After $n$ steps we will have found a path $q=e_{i_n}\dots e_{i_1}$ with $r(q)=w$ and $s(q)=r(e)$. In particular we have
$w\geq s(e)$, and therefore there exists a cycle based at $w$ containing the edge $e$.  Since $e$ is not in $p$ we get
$|CSP(w)|\geq 2$.  Since $w$ is a vertex contained in the cycle $p$, we then get $|CSP(v)|\geq 2$, contrary to the
definition of the set $V_1$. \endproof

\re{\locallymatricialnotpurelyinfinite}{\bf Lemma} Suppose $A$ is a union of finite dimensional subalgebras. Then $A$
is not purely infinite.  In fact, $A$ contains no infinite idempotents.
\endre
\proof  It suffices to show the second statement. So just suppose $e = e^2 \in A$ is infinite.  Then $eA$ contains a
proper direct summand isomorphic to $eA$, which in turn, by definition and a standard argument, is equivalent to the
existence of elements $g,h,x,y \in A$ such that $g^2 = g, h^2 = h, gh = hg = 0, e = g+h, h\neq 0, x \in eAg, y\in gAe$
with $xy = e$ and $yx = g$.  But by hypothesis the five elements $e,g,h,x,y$ are contained in a finite dimensional
subalgebra $B$ of $A$, which would yield that $B$ contains an infinite idempotent, and thus contains a non-artinian
right ideal, which is impossible.
\endproof

\re{\cornernotpurelyinfinite}{\bf Proposition} Let $E$ be a graph. Suppose that
$w\in E^0$ has the property that, for every $v\in E^0$, $w\leq v$
implies $v\in V_0$. Then the corner algebra $wL(E)w$ is not purely infinite.
\endre

\proof Consider the graph $H=(H^0,H^1,r,s)$ defined by $H^0:=\{v:w\leq v\}$, $H^1:=s^{-1}(H^0)$, and $r,s$ induced by
$E$. The only nontrivial part of showing that $H$ is a well defined graph is verifying that $r(s^{-1}(H^0))\subseteq
H^0$. Take $z\in H^0$ and $e\in E^1$ such that $s(e)=z$. But we have $w\leq z$ and thus $w\leq r(e)$ as well, that is,
$r(e)\in H^0$.

Using that $H$ is acyclic, along with the same argument as given in Theorem \ultramatricial, we have that $L(H)$ is a
subalgebra of $L(E)$.  Thus Proposition \ultramatricial\ applies, which yields that $L(H)$ is the union of finite
dimensional subalgebras, and therefore contains no infinite idempotents by Lemma \locallymatricialnotpurelyinfinite.
As $wL(H)w$ is a subalgebra of $L(H)$, it too contains no infinite idempotents, and thus is not purely infinite.

 We claim that $wL(H)w=wL(E)w$.  To see this, given $\alpha
=\sum p_iq_i^*\in L(E)$, then $w\alpha w=\sum p_ {i_j}q_{i_j}^*$ with $s(p_{i_j})=w=s(q_{i_j})$ and therefore
$p_{i_j},q_{i_j}\in L(H)$.  Thus $wL(E)w$ is not purely infinite as desired.
\endproof

We thank P. Ara for indicating the following result, which will provide the direction of proof for our main theorem.  A
right $A$-module $T$ is called {\bfm directly infinite} in case $T$ contains a proper direct summand $T'$ such that
$T'\cong T$. (In particular, the idempotent $e$ is infinite precisely when $eA$ is directly infinite.) Recall that a
ring $A$ has {\bfm local units} if for every finite subset $\{x_1,\dots,x_n\}\subseteq A$ there exists $e=e^2\in A$
with $x_i\in eAe$ for every $i=1,\dots,n$.

\re{\PereAraLemma}{\bf Proposition}  Let $A$ be a ring with local units.  The following are equivalent:

\item{(i)} $A$ is purely infinite simple.
\item{(ii)} $A$ is simple, and for each nonzero finitely generated projective right $A$-module $P$, every nonzero submodule $C$ of
$P$ contains a direct summand $T$ of $P$ for which $T$ is directly infinite.  (In particular, the property `purely
infinite simple' is a Morita invariant of the ring.)
\item{(iii)} $wAw$ is purely infinite simple for every nonzero idempotent $w\in A$.
\item{(iv)} $A$ is simple, and there exists a nonzero idempotent $w$ in $A$ for which $wAw$ is purely infinite simple.
\item{(v)} $A$ is not a division ring, and $A$ has the property that for every pair of nonzero elements $\alpha, \beta$
in $A$ there exist elements $a,b$ in $A$ such that $a\alpha b = \beta$.
\endre

\proof (i) $\Leftrightarrow$ (ii).  Suppose $A$ is purely infinite simple.  Let $P$ be any nonzero finitely generated
projective right $A$-module. Then $P$ is a generator for $Mod-A$, as follows.  Since $A$ generates $Mod-A$ and $P$ is
finitely generated we have an integer $n$ such that $P\oplus P'\cong A^n$ as right $A$-modules.  Again using that $P$
is finitely generated, and using that $A$ has local units, we have that $P$ is isomorphic to a direct summand of a
right $A$-module of the form $f_1A \oplus ... \oplus f_tA$, where each $f_i$ is idempotent.  But this gives
$Hom_A(P,f_1A \oplus \dots \oplus f_tA)\neq 0$, which in turn gives $0 \neq Hom_A(P,A^t) \cong (Hom_A(P,A))^t$, so that
$Hom_A(P,A) \neq 0$. But $\Sigma \{a\in A\mid a=g(p)$ for some $p \in P$ and some $g\in Hom_A(P,A)\}$ is then a nonzero
two-sided ideal of $A$, which necessarily equals $A$ as $A$ is simple.   Now let $e = e^2\in A$.  Then $e =
\sum_{i=1}^r g_i(p_i)$ for some $p_i \in P$ and $g_i\in Hom_A(P,A)$, which gives that $\lambda_e \circ \oplus g_i:
P^{r} \rightarrow A \rightarrow eA$ is a surjection.  Since $P$ generates $eA$ for each idempotent $e$ of $A$, we
conclude that $P$ generates $Mod-A$.

This observation allows us to argue exactly as in the proof of \cite{\AGP, Lemma 1.4 and Proposition 1.5} that if
$e=e^2\in A$, then there exists a right $A$-module $Q$ for which $eA\cong P\oplus Q$. Since $A$ is purely infinite,
there exists  an infinite idempotent $e\in A$.  The indicated isomorphism yields that any submodule $C$ of $P$ is
isomorphic to a submodule $C'$ of $eA$, so that by the hypothesis that $A$ is purely infinite we have that $C'$
contains a submodule $T'$ which is directly infinite, and for which $T'$ is a direct summand of $eA$. But by a standard
argument, any direct summand of $eA$ is equal to $fA$ for some idempotent $f \in A$, so that $T' = fA$ for some
infinite idempotent $f$ of $A$.  Let $T$ be the preimage of $T'$ in $P\oplus Q$ under the isomorphism. Then $T$ is
directly infinite, and since $fA$ is a direct summand of $eA$ we have that $T$ is a direct summand of $P\oplus Q$ which
is contained in $P$, and hence $T$ is a direct summand of $P$.

By \cite{\AM, Proposition 3.3}, the lattice of two-sided ideals of Morita equivalent rings are isomorphic, so that any
ring Morita equivalent to a simple ring is simple. Therefore, since the indicated property is clearly preserved by
equivalence functors, we have that `purely infinite simple' is a Morita invariant.

For the converse, let $I$ be a nonzero right ideal of $A$. We show that $I$ contains an infinite idempotent.   Let
$0\neq x \in I$, so that $xA\leq I$.   But $x=ex$ for some $e = e^2 \in A$, so $xA\leq eA$. So by hypothesis,   $xA$
contains a nonzero direct summand $T$ of $eA$, where $T$ is directly infinite.   But as noted above we have that $T =
fA$ for $f=f^2 \in A$, where $f$ is infinite.  Thus $f\in T \leq xA \leq I$ and we are done.

(ii) $\Rightarrow$ (iii).  Since we have established the equivalence of (i) and (ii), we may assume $A$ is purely
infinite simple.   Then the simplicity of $A$ gives that $AwA=A$ for any nonzero idempotent $w\in A$, which yields by
\cite{\AM, Proposition 3.5} that $A$ and $wAw$ are Morita equivalent, so that (iii) follows immediately from (ii).

(iii) $\Rightarrow$ (iv). It is tedious but straightforward to show that if $A$ is any ring with local units, and $wAw$ is a simple (unital) ring for every nonzero idempotent $w$ of $A$, then $A$ is simple.

(iv) $\Rightarrow$ (i).   Since $A$ is simple we get $AwA=A$, so that $A$ and $wAw$ are Morita equivalent by the
previously cited \cite{\AM, Proposition 3.5}.

Thus we have established the equivalence of statements (i) through (iv).

(i) $\Rightarrow$ (v).   Suppose $A$ is purely infinite simple.  Then $A$ is not left artinian, so that $A$ cannot be a
division ring. Now choose nonzero $\alpha, \beta \in A$.   Then there exists a nonzero idempotent $w \in A$ such that
$\alpha, \beta \in wAw$. But $wAw$ is purely infinite simple by (i) $\Leftrightarrow$ (iii), so by \cite{\AGP, Theorem
1.6} there exist $a',b' \in wAw$ such that $a' \alpha b' = w$.  But then for $a = a', b = b' \beta$ we have $a \alpha b
= \beta$. Conversely, suppose $A$ is not a division ring, and that $A$ satisfies the indicated property.  Since $A$ is
not a division ring and $A$ is a ring with local units, there exists a nonzero idempotent $w$ of $A$ for which $wAw$ is
not a division ring.  Let $\alpha \in wAw$.  Then by hypothesis there exist $a',b'$ in $A$ with $a' \alpha b' = w$. But
since $\alpha \in wAw$, by defining $a = wa'w$ and $b = wb'w$ we have $a \alpha b = w$.  Thus another application of
\cite{\AGP, Theorem 1.6} (noting that $w$ is the identity of $wAw$) gives the desired conclusion.

(v) $\Rightarrow$ (iv).  The indicated multiplicative property yields that any nonzero ideal of
$A$ will contain a set of local units for $A$, so that $A$ is simple.
 Since $A$ is not a division ring and $A$ has local units there exists a nonzero idempotent
$w$ of $A$ such that $wAw$ is not a division ring.  Let $\alpha, \beta \in wAw$; in particular, $w\alpha w = \alpha$
and $w\beta w = \beta$.  By hypothesis there exists $a,b \in A$ such that $a\alpha b = \beta$.  But then
$(waw)\alpha(wbw) = w\beta w = \beta$, which yields that $wAw$ is purely infinite simple by \cite{\AGP, Theorem 1.6}.
\endproof

We now have all the necessary ingredients in hand to prove the main result of this
article.

\re{\purelyinfinite}{\bf Theorem} Let $E$ be a graph. Then $L(E)$ is purely infinite simple if and only if $E$ has the
following properties.

\item{(i)} The only hereditary and saturated subsets of $E^0$ are $\emptyset$ and $E^0$.
\item{(ii)} Every cycle in $E$ has an exit.
\item{(iii)} Every vertex connects to a cycle.

\endre
\proof  First, assume (i), (ii) and (iii) hold.  By Theorem \simple\ we have that $L(E)$ is simple. By Proposition
\PereAraLemma\ it suffices to show that $L(E)$ is not a division ring, and that for every pair of elements $\alpha,
\beta$ in $L(E)$ there exist elements $a,b$ in $L(E)$ such that $a\alpha b = \beta$. Conditions (ii) and (iii) easily
imply that $|E^1|>1$, so that $L(E)$ has zero divisors, and thus is not a division ring.

We now apply Proposition \reducetovertex\ to find $\overline{a}$, $\overline{b} \in L(E)$ such that $\overline{a}\alpha
\overline{b}=w\in E^0$. By condition (iii), $w$ connects to a vertex $v\not\in V_0$. Either $w=v$ or there exists a
path $p$ such that $r(p)=v$ and $s(p)=w$. By choosing $a' = b' = v$ in the former case, and $a' = p^*, b' = p$ in the
latter, we have produced elements $a', b' \in L(E)$ such that $a'wb' = v$.

An application of Lemma \voneempty\ yields that $v\in V_2$, so there exist $p,q \in CSP(v)$ with $p\neq q$. For any
$m>0$ let $c_m$ denote the closed path $p^{m-1}q$. Using \cite{\AA, Lemma 2.2}, it is not difficult to show that $c^*_m
c_n=\delta_{mn}v$ for every $m,n>0$.

Now consider any vertex $v_l \in E^0$. Since $L(E)$ is simple, there exist
$\{a_i,b_i\in L(E)\mid 1\leq i \leq t\}$
such that $v_l=\sum_{i=1}^t a_i v b_i$. But by defining  $a_l=\sum_{i=1}^t a_i
c^*_i$ and $b_l=\sum_{j=1}^t c_j b_j$,
we get $$a_l v b_l =\left(\sum_{i=1}^t a_i c^*_i\right)v\left(\sum_{j=1}^t c_j
b_j\right)=\sum_{i=1}^t a_i c^*_i v c_i
b_i=v_l.$$

Now let $s$ be a left local unit for $\beta$ (i.e., $s \beta = \beta$), and write
$s=\sum_{v_l\in S} v_l$ for some finite subset of vertices $S$.
By letting $\widetilde{a}=\sum_{v_l\in S} a_l c^*_l$ and $\widetilde{b}=\sum_{v_l\in
S} c_l b_l$, we get
$$\widetilde{a} v \widetilde{b}= \sum_{v_l\in S} a_l c^*_l v c_l b_l =\sum_{v_l\in
S} v_l=s.$$

Finally, letting $a=\widetilde{a} a' \overline{a}$ and $b= \overline{b} b'
\widetilde{b}\beta$, we have that $a \alpha b =\beta$ as desired.

For the converse, suppose that $L(E)$ is purely infinite simple. By Theorem \simple\ we have (i) and (ii).
 If (iii) does not hold, then there exists a vertex $w\in E^0$ such
that $w\leq v$ implies $v\in V_0$. Applying Proposition \cornernotpurelyinfinite\ we get that $wL(E)w$ is not purely
infinite.  But then Proposition \PereAraLemma\ implies that $L(E)$ is not purely infinite, contrary to hypothesis.
\endproof

\re{\Leavittalgebra}{\bf Examples}

{\rm \item{(i)} Let $E$ be the graph defined in Example \Leavittalgebraexamples\ (i). Then $L(E) \cong M_n(K)$ which of
course is simple, but not purely infinite since no vertex in $E^0$ connects to a cycle.
\item{(ii)} Let $n\geq 2$. Let  $E$ be the graph defined in Example \Leavittalgebraexamples\ (ii). Then $L(E)\cong L(1,n)$, the Leavitt algebra.
Since $n\geq 2$ we see that all the hypotheses of Theorem \purelyinfinite\ are satisfied, so that $L(1,n)$ is purely
infinite simple.
\item{(iii)} Let $E$ be the graph having $E^0 = \{v,w\}$ and $E^1=\{e,f,g\}$, where $s(e)=s(f)=v$, $r(e)=r(f)=w$, $s(g)=w, r(g)=v$.
Then  $E$ satisfies the hypotheses of Theorem \purelyinfinite, so that $L(E)$ is purely infinite simple.
\endre}
\endre

Let $L(1,n)$ denote the Leavitt algebra described in Example \Leavittalgebraexamples\ (ii). We complete this article by
providing a realization of the purely infinite simple algebra $M_m(L(1,n))$ as a Leavitt path algebra $L(E)$ for a
specific graph $E$.

\re{\example}{\bf Proposition} Let $n\geq 2$ and $m\geq 1$. We define the graph $E_n^m$ by setting
$E^0:=\{v_1,\dots,v_m\}$, $E^1:=\{f_1,\dots,f_n,e_1,\dots,e_{m-1}\}$, $r(f_i)=s(f_i)=v_m$ for $1\leq i \leq n$,
$r(e_i)=v_{i+1}$, and $s(e_i)=v_i$ for $1\leq i \leq m-1$. Then $L(E_n^m)\cong M_m(L(1,n))$.
\endre

\proof We define $\Phi:K[E^0\cup E^1 \cup (E^1)^*]\to M_m(L(1,n))$ on the generators
by $$\eqalign{
 \Phi(v_i) & = e_{ii} \hbox{ for }1\leq i \leq m  \cr
 \Phi(e_i) & = e_{ii+1} \hbox{ and }\Phi(e_i^*)=e_{i+1i} \hbox{ for }1\leq i \leq
m-1  \cr
 \Phi(f_i) & = y_ie_{mm} \hbox{ and }\Phi(f_i^*)=x_ie_{mm} \hbox{ for }1\leq i \leq
n \cr}$$ and extend linearly and multiplicatively to obtain a $K$-homomorphism. We now verify that $\Phi$ factors
through the ideal of relations in $L(E_n^m)$.

First, $\Phi(v_iv_j-\delta_{ij}v_i)=e_{ii}e_{jj}-\delta_{ij}e_{ii}=0$. If we consider the relations $e_i-e_ir(e_i)$
then we have $\Phi(e_i-e_ir(e_i))=\Phi(e_i-e_iv_{i+1})=e_{ii+1}-e_{ii+1}e_{i+1i+1}=0$, and analogously
$\Phi(e_i-s(e_i)e_i)=0$. For the relations $f_i-f_ir(f_i)$ we get
$\Phi(f_i-f_ir(f_i))=\Phi(f_i-f_iv_m)=y_ie_{mm}-y_ie_{mm}e_{mm}=0$, and similarly $\Phi(f_i-s(f_i)f_i)=0$. With similar
computations it is easy to also see that
$\Phi(e^*_i-e^*_ir(e^*_i))=\Phi(e^*_i-s(e^*_i)e^*_i)=\Phi(f^*_i-f^*_ir(f^*_i))=\Phi(f^*_i-s(f^*_i)f^*_i)=0$.

We now check the Cuntz-Krieger relations. First,
$\Phi(e^*_ie_j-\delta_{ij}r(e_j))=\Phi(e^*_ie_j-\delta_{ij}v_{j+1})=e_{i+1i}e_{jj+1}-\delta_{ij}e_{j+1j+1}=
 \delta_{ij}e_{i+1j+1}-\delta_{ij}e_{j+1j+1}=0$. Second,
 $\Phi(f^*_if_j-\delta_{ij}r(f_j))=\Phi(f^*_if_j-\delta_{ij}v_m)=x_ie_{mm}y_je_{mm}-\delta_{ij}e_{mm}=0$,
because of the relation (1) in $L(1,n)$. Finally,
 $\Phi(f^*_ie_j-\delta_{f_i,e_j}r(e_j))=\Phi(f^*_ie_j-0v_{j+1})=\Phi(f^*_ie_j)=x_ie_{mm}e_{jj+1}=0$,
and similarly
 $\Phi(e^*_if_j-\delta_{e_i,f_j}r(f_j))=0$.

With CK2 we have two cases. First, for $i<m$, $\Phi(v_i-e_ie^*_i)=e_{ii}-e_{ii+1}e_{i+1i}=0$. And for $v_m$ we have
$\Phi(v_m-\sum_{i=1}^n f_if^*_i)=e_{mm}-\sum_{i=1}^n y_ie_{mm}x_ie_{mm}=0$, because of the relation (2) in $L(1,n)$.

This shows that we can factor $\Phi$ to obtain a $K$-homomorphism of algebras $\Phi:L(E_n^m)\to M_m(L(1,n))$. We will
see that $\Phi$ is onto. Consider any matrix unit $e_{ij}$ and $x_k\in L(1,n)$. If we take the path $p=e_i\dots
e_{n-1}f^*_ke^*_{n-1}\dots e^*_j\in L(E_n^m)$ then we get $\Phi(p)=e_{ii+1}\dots e_{n-1n}(x_ke_{nn})e_{nn-1}\dots
e_{j+1j}=x_ke_{ij}$. Similarly $\Phi(e_i\dots e_{n-1}f_ke^*_{n-1}\dots e^*_j)=y_ke_{ij}$. In this way we get that all
the generators of $M_m(L(1,n))$ are in $Im(\Phi)$.

Finally, using the same ideas as those presented in \cite{\AA, Corollary 3.13 (i)}, we see that $E_n^m$
satisfies the conditions of Theorem
\simple, which yields the simplicity of $L(E_n^m)$. This implies that $\Phi$ is necessarily
injective, and therefore an
isomorphism.
 \endproof

\vskip 10pt plus 3 pt\noindent{\bfg Acknowledgments}\par\nobreak

The authors are grateful to P. Ara and E. Pardo for many valuable correspondences. The second author was partially
supported by the MCYT and Fondos FEDER, BFM2001-1938-C02-01, the ``Plan Andaluz de Investigaci\'on y Desarrollo
Tecnol\'ogico", FQM 336 and by a FPU fellowship by the MEC (AP2001-1368). This work was done while the second author
was a Research Scholar at the University of Colorado at Colorado Springs supported by a ``Estancias breves" FPU grant.
The second author thanks this host center for its warm hospitality.

\medskip





\vskip 10pt plus 3 pt\noindent{\bfg References}\par\nobreak

[\AA] {\rms G. Abrams, G. Aranda Pino, The Leavitt path algebra of a graph, J. Algebra (to appear).}

[\AM]  {\rms P.N. \'Anh, L. M\'arki, Morita equivalence for rings without identity, Tsukuba J. Math 11(1) (1987) 1-16.}

[\Ados] {\rms P. Ara, The exchange property for purely infinite simple rings, Proc. A.M.S. 132 (2004) 2543-2547.}

[\AGGP] {\rms P. Ara, M.A. Gonz\'alez-Barroso, K.R. Goodearl, E. Pardo, Fractional skew monoid rings, J. Algebra 278
(2004) 104-126.}

[\AGP] {\rms P. Ara, K.R. Goodearl, E. Pardo, $K_0$ of purely infinite simple regular rings, K-Theory 26 (2002)
69-100.}

[\AMP] {\rms P. Ara, M.A. Moreno, E. Pardo, Nonstable K-theory for graph algebras, Submitted for publication.}

[\BPRS] {\rms T. Bates, D. Pask, I. Raeburn, W. Szyma\'nski, The C*-algebras of row-finite graphs, New York J. Math. 6
(2000) 307-324.}

[\KPR] {\rms A. Kumjian, D. Pask, I. Raeburn, Cuntz-Krieger algebras of directed graphs, Pacific J. Math. 184 (1)
(1998) 161-174.}

[\Lone] {\rms W.G. Leavitt, The module type of a ring, Trans. A.M.S. 42 (1962) 113-130.}

[\Ltwo] {\rms W.G. Leavitt, The module type of homomorphic images, Duke Math. J. 32 (1965) 305-311.}

[\P]  {\rms N.C. Phillips, A classification theorem for nuclear purely infinite simple C*-algebras, Doc. Math. 5 (2000)
49-114.}

\end